\documentclass[12pt]{amsart}

\usepackage{amsmath}

\usepackage{amssymb}

\usepackage{latexsym}

\usepackage{amsthm}

\begin{document}

\begin{center}
\Large{\bf Further results on some  singular linear stochastic
differential equations}
\end{center}

\vspace{1cm}

\begin{center}
{\bf Larbi Alili$^{(1)}$} \hspace{5mm} {\bf Ching-Tang Wu$^{(2)}$}
\end{center}

\vspace{5mm}

\begin{center}
{\small {\sc ABSTRACT}} \\
\vspace{4mm}
\begin{minipage}{14cm}
A class of Volterra transforms, preserving the Wiener measure,
with kernels of Goursat type is considered. Such kernels satisfy a
self-reproduction property. We provide some results on the
inverses of the associated Gramian matrices which lead to a new
self-reproduction property. A connection to the classical
reproduction property is given. Results are then applied to the
study of a class of singular linear stochastic differential
equations together with the corresponding decompositions of
filtrations. The studied equations are viewed as non-canonical
decompositions of some generalized bridges.
\end{minipage}
\end{center}

\vspace{1cm}


\vspace{1cm}


\theoremstyle{plain}
   \newtheorem{Theorem}{Theorem}[section]
   \newtheorem{Corollary}{Corollary}[section]
   \newtheorem{Lemma}{Lemma}[section]
   \newtheorem{Proposition}{Proposition}[section]

\theoremstyle{definition}
   \newtheorem{Definition}{Definition}[section]
   \newtheorem{Example}{Example}[section]
   \newtheorem{Remark}{Remark}[section]

\def\real{\mathbb{R}}
\def\rational{\mathbb{N}}
\def\measure{\mathbb{P}}
\def\preal{{\mathbb{R}}^+}
\def\ii{{I\hspace{-0.5mm}I}}
\def\iii{{I\hspace{-0.5mm}I\hspace{-0.5mm}I}}
\def\sign{\mbox{\rm{sign}}}
\newcommand\spn{\mbox{\rm{span}}}
\def\dim{\mbox{\rm{dim}}}
\def\span{\mbox{\rm{span}}}
\def\dim{\mbox{\rm{dim}}}


\section{Introduction and preliminaries}
Gaussian enlargement of filtrations has been extensively studied
between the late 70's and the early 90's, see \cite{deheuvels},
\cite{jeulin80}, \cite{jeulinyor85}, \cite{jeulin88} and the
references therein. Results stemming from the Gaussian nature of
the underlying generalized Gaussian bridges are of interest not
only in probability, also in financial mathematics, since they
have appeared in an insider trading model developed in
\cite{jurgen98a} and \cite{karatzas96}. Transforms of Volterra
type allow to construct interesting families of Gaussian
processes. Volterra-transforms are classified, both from the
theory and applications points of view, according to whether their
kernels are square-integrable or not. Those with square-integrable
kernels play a crucial role in the study of equivalent Gaussian
measures, stochastic linear differential equations and the linear
Kalman-Bucy filter, see \cite{hida93} and \cite{kallianpur80}. To
our knowledge, comparably, less interest was given to Volterra
transforms with non-square-integrable kernels. Such transforms
naturally appear, for instance,  in non-canonical representations
of some Gaussian processes. They also appear if one forces such
transforms to preserve the Wiener measure. The most known examples
have corresponding  kernels of Goursat type. A few nontrivial ones
originate from P. L\'evy, see \cite{Levy-56}, \cite{Levy-57}, and
serve as a standard reference for showing the importance of the
canonical decomposition of semi-martingales. Such constructions
have been enriched by people from the Japanese school, see
\cite{hhm}, \cite{Hida-60} and \cite{Hitsuda-68}.

Let us now fix the mathematical setting and summarize results of
this paper. We take $B:=(B_t, t\geq 0)$ to be a standard Brownian
motion, defined on a complete probability space $(\Omega,
{\mathcal F}, \mathbb{P}_0)$. Denote by $\{\mathcal{F}_t^B, t\geq
0\}$ the filtration it generates. Let $f = (f_1, \cdots, f_n)^*
\in L_{loc}^2(\real_+)=\{ h; \int_0^t h^2(s) \, ds<\infty, \hbox{
for all } t\in [0, \infty)\}$, where $*$ stands for the transpose
operator and $n$ is a natural number. Although  some of our
results extend readily to the cases when $n=\infty$, to simplify
the study, we only consider the cases where $n$ is finite. We
assume that, for any fixed $t>0$, the covariance matrix $m_t$, of
the Gaussian random variable $\int_0^t f^*(s) \, d B_s$, is
invertible, i.e., the Gramian matrix $m_t=\int_0^t f(s) \cdot
f^*(s) \, ds$ has an inverse $\alpha_t$. We emphasize that, under
the aforementioned condition, it is not difficult to see that
$\alpha_t\rightarrow \alpha_{\infty}$, as $t\rightarrow \infty$,
where $\alpha_{\infty}$ is a finite matrix. Furthermore, for any
$i$, $(\alpha_{\infty})_{ij}=0$ for all $j$, if and only if
$\|f_i\|:=(\int_0^{\infty}f_i^2(s) \, ds)^{1/2}=\infty$. With
$\phi(t) = \alpha_t \cdot f(t)$ for $t>0$, we shall establish in
Theorem \ref{Theorem29} that $(\alpha_{t}, t>0)$ is given in terms
of $\phi$ by $\alpha_t= \int_t^{\infty} \phi(u) \cdot \phi^*(u) \,
du + \alpha_{\infty}$, for any $t > 0$. This relation has its own
right of importance in this work and may have interesting
applications to other fields where Gramian matrices together with
their inverses are of prime importance, see for instance
\cite{Berlinet-Thomas} and the references therein. In particular,
we also refer to \cite{Andrews-Askey-Roy} for applications to the
theory of special functions and to \cite{Aronszajn-50} and
\cite{Atteia-92} for applications to reproducing kernel-Hilbert
spaces and spline functions.

We define the Volterra transform $\Sigma$, associated to a
Volterra kernel $k$, on the set of continuous semi-martingales $X$
such that
\begin{equation}\label{condition}
 \lim_{\varepsilon \to
0} \int_{\varepsilon}^t  \int_0^v k(u,v) \, dX_u \, dv < \infty,
\qquad  0<t<\infty\quad \mbox{\rm{a.s.}},
\end{equation}
 by
\begin{equation}\label{Volterra}
\Sigma(X)_{t} = X_{t} - \int_0^{t} \! \int_0^u k(u,v) \, dX_v \,
du, \qquad 0<t<\infty.
\end{equation}
Following \cite{hhm}, the kernel $k(t,s) = \phi^*(t)\cdot f(s)$,
for $0< s \leq t<+\infty$ is a self-reproducing Volterra kernel.
That is equivalent to saying that $\Sigma$, when applied to the
Brownian motion $B$, satisfies the following two conditions:
\begin{itemize}
 \item[{\bf (i)}] $\Sigma(B)$ is a standard Brownian motion;
 \item[{\bf (ii)}] For any fixed $t\geq 0$, ${\mathcal F}_t^{\Sigma(B)}$
 is independent of $\int_0^t f(u) \, dB_u $.
\end{itemize}
Existence of  $\Sigma(B)$  may be justified by using a generalized
Hardy inequality discovered in \cite{hhm}, see Remark \ref{Hardy}
given below.   We call $k$ and $\Sigma$, respectively, a {\it
Goursat-Volterra kernel} and {\it transform}, with reproducing
basis $f$. The dimension of $Span\{f\}$ is called the order of the
Goursat-Volterra kernel $k$. This terminology is formally fixed in
Definition \ref{Definition-g-v-k}.

Next, we bring our focus on conditions {\bf (i)} and {\bf (ii)}
and think of them in terms of enlargement of filtrations and
stochastic differential equations. Condition {\bf (ii)} says that
the orthogonal decomposition
\begin{equation}
\label{eqn:C68} {\mathcal F}_t^B = {\mathcal F}_t^{\Sigma(B)}
\otimes \sigma \left( \int_0^t f(u) \, dB_u \right)
\end{equation}
holds true, for any $t \geq 0$. Here,  by $\mathcal{F} \otimes \,
\mathcal{G}$ we mean $\mathcal{F} \vee \mathcal{G}$ with
independence between $\mathcal{F}$ and $\mathcal{G}$.  We shall
show that, for Goursat-Volterra transforms, equation
(\ref{eqn:C68}) can in fact be rewritten as
\begin{equation}
\label{eqn:C70} {\mathcal F}_t^B = {\mathcal F}_t^{\Sigma(B)}
\otimes \sigma \left( Y - \int_t^{\infty} \phi(u) \,
d\Sigma(B)_u\right)
\end{equation}
valid for any $t \geq 0$, where $Y=(Y_1, \cdots, Y_N)^*$ is a
Gaussian random vector which is independent of ${\mathcal
F}_{\infty}^{\Sigma(B)}$ with covariance matrix $E[Y \cdot Y^*] =
\alpha_{\infty} = \lim_{t \to \infty} \alpha_t$ in case
$\alpha_{\infty} \not\equiv 0$, and $Y \equiv 0$ otherwise. We
allow here $Y$ to have some null or constant components. Going
back to condition {\bf (i)}, we observe that the determination of
all continuous semi-martingales which satisfy it amounts to
solving equation
\begin{equation}\label{ssde}
X_t = W_t + \int_0^{t} \! \int_0^s \phi^*(s)\cdot f(u) \, dX_u \,
ds, \, \quad X_0 = 0, \qquad t>0,
\end{equation}
considered on a possibly enlarged probability space, where $W$ is
a standard Brownian motion. Note that we only assume
\begin{equation}\label{condition-2}
\lim_{\varepsilon \to 0} \int_{\varepsilon}^t \! \int_0^v
\phi^*(v) \cdot f(u) \, dX_u \, dv < \infty,  \qquad  0<t<\infty
\quad \mbox{\rm{a.s.}},
\end{equation}
and the latter is not absolutely convergent. Because of the
singularity at time $0$, we call (\ref{ssde}) a singular linear
stochastic differential equation.  If we take $W=\Sigma(B)$ then,
by construction, the original Brownian motion $B$ is one solution.
A second one coincides with the associated $f$-generalized bridge
on the interval of its finite life-time, introduced in
\cite{alili00}. It follows that the Goursat-Volterra transform
$\Sigma$, when defined as above, is not invertible in the sense
that (\ref{ssde}) has many solutions. This is not a surprising
fact. Indeed, $k$ being a self-reproducing kernel implies that it
is not square-integrable, as seen in \cite{fwy99}. Next, Theorem
\ref{thm25} deals with the investigation of all continuous
semi-martingale solutions to (\ref{ssde}). In particular, we show
that a necessary and sufficient condition for the existence of a
strong solution that is Brownian and ${\mathcal
F}_{\infty}^B$-measurable is $\alpha_{\infty}\equiv 0$. In that
case ${\mathcal F}_{\infty}^{\Sigma(B)} = {\mathcal
F}_{\infty}^B$. When $\alpha_{\infty}\not\equiv 0$, Theorem
\ref{thm25} concludes that there exists still a strong solution
which is a Brownian motion, in an enlarged space, that involves an
independent centered Gaussian vector $Y$ with covariance matrix
$\alpha_{\infty}$. Another natural question is a characterization
of all continuous semi-martingales that satisfy both conditions
{\bf (i)} and {\bf (ii)}. This is partially solved in Theorem
\ref{thm-harmonic} for the case $\alpha_{\infty}\equiv 0$ and the
analysis exhibits some connections to certain space-time harmonic
functions. The latter are functions $h \in C^{1,2}\left( \real_+
\times \real^n, \real_+\right)$ such that $ h(\cdot, \int_0^{.}
f^*(s) \, dB_s)$ is a continuous $(\mathbb{P}_0,{\mathcal
F})$-martingale with expectation $1$, where $\mathbb{P}_0$ stands
for the Wiener measure.

The main results of this paper extend a part of the first chapter
of \cite{yor92} and some results found in \cite{jeulin88}. Our
work offers explicit examples of conditionings and conditioned
stochastic differential equations introduced and studied in
\cite{Baudoin}. Furthermore, singular equations of type
(\ref{ssde}) and the progressive enlargement of filtration given
in Corollary \ref{decomposition} can easily be applied to insider
trading models elaborated in \cite{jurgen98a}, \cite{Baudoin} and
\cite{{karatzas96}}.

\section{\bf Goursat-Volterra kernels and transforms }\label{sec35}

To a Brownian motion $B$ we associate  the centered Gaussian
process $\Sigma(B)$ defined by (\ref{Volterra}), which we assume
is well-defined, where $k$ is a continuous Volterra kernel. That
is to say that $k:\mathbb{R}_+^2\rightarrow \mathbb{R}$ satisfies
\begin{equation*}
k(u,v)=0,\qquad 0< u \le v <\infty
\end{equation*}
and is continuous on $\{(u,v) \in (0,+\infty)\times (0,+\infty): u
> v \}$. We know from \cite{fwy99} that $\Sigma$ preserves the
Wiener measure, or $\Sigma(B)$ is a Brownian motion, if and only
if $k$ satisfies the self-reproducing property
\begin{equation}
\label{eqn:K10} k(t,s) = \int_0^s k(t,u) k(s,u) \, du, \qquad 0< s
\leq t<\infty.
\end{equation}
For a connection with reproducing kernels, in the usual sense, we
refer to end of this section. Observe that (\ref{Volterra}), when
applied to $B$, can be viewed as the semi-martingale decomposition
of $\Sigma(B)$ with respect to the filtration $(\mathcal{F}^{B}_t,
t\geq 0)$. Now, as a consequence of the Doob-Meyer decomposition
of $\Sigma(B)$  in its own filtration, we must have the strict
inclusion
\begin{equation*}
{\mathcal F}_t^{\Sigma(B)}
 \subsetneqq {\mathcal F}_t^{B},\qquad  0<t<\infty.
\end{equation*}
It is shown in \cite{jeulin88} that the missing information,
called the reproducing Gaussian space, is given in the orthogonal
decomposition
$$ {\mathcal F}_t^{B} = {\mathcal F}_t^{\Sigma(B)} \otimes
\sigma(\Gamma_t^{(k)}),$$
where
\begin{equation*} \Gamma_t^{(k)} = \left\{ \int_0^t f(u) \, dB_u; f\in L^2\left((0,t]
\right),\: f(s)=\int_0^s k(s,u)f(u) \, du\quad \hbox{ a.e.}
\right\}
\end{equation*}
for any $t>0$. Given a kernel $k$, it is not an easy task to
determine a basis of $\Gamma_t^{(k)}$ for each fixed $t>0$,
because this amounts to solving explicitly the integral equation
\begin{equation*}
 f(t)
= \int_0^s k(t,u) f(u) \, du,\qquad 0<t<\infty.
\end{equation*}
It is easier to fix the family of spaces $(\Gamma_t^{(k)}, t>0)$
and work out the corresponding Volterra kernel. This procedure, in
fact, corresponds to desintegrating the Wiener measure over the
interval $[0,t]$, for any fixed $t>0$, along $\Gamma_t^{(k)}$.
Recall that a Goursat kernel is a kernel of the form
\begin{equation*}
k(t,s) =\phi^*(t) \cdot f(s),\qquad 0< s\leq t<\infty,
\end{equation*}
 where $\phi=(\phi_1, \cdots,
\phi_n)^*$ and $f=(f_1, \cdots, f_n)^*$ are two vectors of
functions defined on $(0,\infty)$ and $n \in \mathbb{N}$. For such
kernels it is natural to introduce the following definition.

\begin{Definition}\label{Definition-g-v-k}
A Goursat-Volterra transform $\Sigma$ of order $(n_t, t>0)$ is a
Volterra transform preserving the Wiener measure such that, for any
Brownian motion $B$ and  $t
> 0$,  ${\mathcal F}_t^{\Sigma(B)}$ is independent of $\int_0^t
f(u) \, dB_u $ for some vector $f\equiv(f_1, \cdots, f_{n_t})^*$
of $n_t$ linearly independent $L^2_{loc}(\real_+)$ functions. The
associated kernel is called a Goursat-Volterra kernel. The objects
$f$, $Span\{f\}$ and $Span\{\int_0^{\cdot} f(s) \, dB_s\}$ are
called reproducing basis, space and Gaussian space, respectively.
\end{Definition}

Because for each fixed $t>0$, $m_t$ is positive definite, it can
be seen that $t\rightarrow n_t$ is nondecreasing. However, in our
setting, we always take the order to be constant and finite. The
simplest known example of Goursat-Volterra kernels is
$k_1(t,s)=t^{-1}$ and this gives
\begin{equation*}
\Sigma(B)_{.}=B_.-\int_0^{\cdot} \frac{B_u}{u} \, du.
\end{equation*}
 That corresponds
to setting $n=1$ and taking $f_1\equiv 1$. It is observed in
\cite{jeulin88} that $\Sigma$ when iterated takes a remarkably
simple form. That is with $\Sigma^{(0)}=Id$, $\Sigma^{(1)}=\Sigma$
and $\Sigma^{(m)} = \Sigma^{(m-1)} \circ \Sigma$, for $m \geq 2$,
where $\circ$ stands for the composition operation, we have
\begin{equation*}
\Sigma^{(n)}(B)_{\cdot}=\int_0^{.}L_n(\log{\frac{\cdot}{s}}) \,
dB_s,
\end{equation*}
 where $(L_n, n\in\mathbb{N})$ is the sequence of Laguerre
polynomials. As a generalization of the above kernel, we quote the
following result from \cite{hhm}.

\begin{Theorem}[Hibino-Hitsuda-Muraoka, \cite{hhm}]\label{hhm-thm}
Let $f$  be a  vector of $n$ functions of $L^2_{loc}(\real_+)$
such that for any $t>0$ the Gramian matrix $m_{t}=\int_0^{t} f(s)
\cdot f^*(s) \, ds$ has an inverse denoted by $\alpha_t$. Then,
with $\phi(\cdot)=\alpha_{\cdot} \cdot f(\cdot)$, the kernel $k$,
defined by $k(t,s)=0$ if $s>t$ and $k(t,s)=\phi^*(t)\cdot f(s)$
otherwise, is a Goursat-Volterra kernel of order $n$.
\end{Theorem}

For a proof of this result, we refer to \cite{hhm}. Some arguments
of the proof are sketched in Remark \ref{Hardy} given below. In
the remainder of this paper, unless otherwise specified, we work
under the setting of Theorem \ref{hhm-thm}. The objective of the
next result is to obtain an expression of $\alpha_{\cdot}$ in
terms of $\phi(\cdot)$. As a straightforward application, we shall
show that it allows to obtain a new self-reproducing property
satisfied by the kernel $k$. To our knowledge the following result
is not known.

\begin{Theorem} \label{Theorem29} $\alpha_t$ converges to a finite
matrix $\alpha_{\infty}$ as $t\rightarrow \infty$. Moreover, we
have
\begin{equation}\label{expression-alpha}
\alpha_t= \int_t^{\infty}\phi(u)\cdot \phi^*(u) \, du +
\alpha_{\infty}, \qquad 0<t<\infty.
\end{equation}
Consequently, the self-reproduction property
\begin{equation}
 k(t,s) = \int_t^{\infty}
k(u,t) k(u,s) \, du + f^*(t)\cdot \alpha_{\infty}\cdot f(s), \quad
0<s\leq t<\infty,
\end{equation}
holds true.
\end{Theorem}

\begin{proof} Fix $t>0$. Observe that the matrices $\alpha_t$ and $m_t$ are
symmetric positive definite with absolutely continuous entries.
Next, the identity $\alpha_t\cdot m_t=Id_n=m_t\cdot \alpha_t$,
when differentiated, yields  $\alpha'_t\cdot m_t=-\alpha_t\cdot
m'_t$. It follows that
\begin{equation*}
\phi(t)\cdot f^*(t)=\alpha_t\cdot f(t)\cdot f^*(t) =\alpha_t\cdot
m'_t =-\alpha'_t\cdot m_t.
\end{equation*}
Consequently, we have $\alpha'_t=-\phi(t)\cdot f^*(t)\cdot
\alpha_t=-\phi(t)\cdot\phi^*(t)$. For any $1\leq j\leq n$,
$(\alpha'_t)_{j,j}=-\phi_j^2(t)$ is negative. Hence,
$(\alpha_t)_{j,j}$ is decreasing. Because $(\alpha_t)_{j,j}>0$ we
get that $ \int_r^{\infty}\phi_{j}^{2}(s) \, ds <\infty$, $r>0$.
Since, for  $t\geq r$, we can write $\alpha_t =
\alpha_r-\int_r^t\phi(s) \cdot \phi^*(s) \, ds$, by letting $t
\rightarrow +\infty$, we find $\lim_{t\rightarrow \infty}\alpha_t=
\alpha_r-\int_r^{\infty}\phi(s)\cdot \phi^*(s) \,
ds=\alpha_{\infty}$. Thus, $\alpha_{\infty}$ is a matrix with
finite entries. The last statement follows from $k(t,s) = f^*(t)
\cdot \alpha_t \cdot f(s)$ where we use the expression of
$\alpha_t$ given in (\ref{expression-alpha}).
\end{proof}

Self-reproducing kernels, in particular Goursat-Volterra kernels,
are different from but related to kernel systems and reproducing
kernel Hilbert spaces. Our next objective is to outline this
connection. For, let us start by fixing a time interval $[0,t]$,
for some $t>0$. Let the vector $q_t(u):=(q_{m,t}(u), 0<u\leq t;
1\leq m\leq n)$ be formed by the orthonormal sequence associated
to  $f_1, f_2,\cdots, f_n$ over the interval $[0,t]$. This system
is uniquely characterized by
\begin{equation*}
\int_0^t q_{m, t}(r)q_{k, t}(r) \, dr=\delta_{m,k}, \qquad 1\leq
m,k\leq n,
\end{equation*}
with the requirement that for each integer $1\leq m\leq n$,
$q_{m,t}$ is a linear combination of $f_1$, ..., $f_{m}$ with a
positive leading coefficient associated to $f_{m}$. We refer to
Lemma 6.3.1, p. 294, in \cite{Andrews-Askey-Roy} for an expression
of the latter in terms of a determinant. The classical kernel
system is then given by the symmetric kernel
\begin{equation*}
\kappa_{t}(u,v)=q_{t}(u) \cdot q_{t}^{*}(v), \qquad 0<u,v\leq t.
\end{equation*}
 This is a reproducing kernel in the sense that
$$\kappa_{t}(u,v)=\int_0^{t} \kappa_{t}(u,r)\kappa_{t}(v,
r) \, dr,\qquad 0< u, v\leq t.
$$
For $1\leq i,j\leq n$, $(\alpha_t)_{i,j}$ is seen to be the
coefficient of $f_i(u)f_j(v)$ in the expansion of $\kappa_{t}$. To
be more precise, $(\alpha_t)_{i,j}= (b_t \cdot b^*_t )_{i,j}$
where $b$ is an upper diagonal matrix whose entry $(b_t)_{i,k}$ is
the coefficient of $f_i(u)$ in $q_{k,t}(u)$ for $i\leq k$. We
clearly have $\phi_i^2(t)=-2 (b_t' \cdot b^*_t)_{i,i}$ for all $i$
and it would be interesting to express the matrix $b_t$ in terms
of $\phi(t)$. Now, we are ready to state the following result
which proof is omitted.

\begin{Proposition} For each fixed $t>0$, the kernel system
associated to $f$, over the time interval $[0,t]$, is given by
$\kappa_{t}(u,v)=\int_t^{\infty} k(r,u)k(r,v) \, dr +f^*(u) \cdot
\alpha_{\infty} \cdot f(v) $ for $0< u,v \leq t$. In particular,
we have $k(t,s)=\kappa_{t}(t,s)$ for all $0 < s \leq t<\infty$.
\end{Proposition}

\begin{proof} As in the proof of Theorem \ref{Theorem29},
the first part of the result follows from the well-known
relationship $\kappa_{t}(u,v)=f^*(u) \cdot \alpha_t \cdot f(v)$
for any $0 < u, v \leq t$. The second part follows by taking the
limit and using continuity.
\end{proof}

\begin{Remark}
To see an example where $\alpha_{\infty}\not\equiv 0$, let us
discuss the case $n = 2$. Assume that $f_1$ and $f_2$ are two
functions in $L^2_{loc}(\mathbb{R}_+)$. We distinguish four cases
and three different forms for $\alpha_{\infty}$. The first
corresponds to $\alpha_{\infty}\equiv 0$ when
$\|f_1\|=\|f_2\|=+\infty$. The second corresponds to case when $
\|f_1\|$ and $\|f_2\|$ are finite which implies that
$\alpha_{\infty}$ is positive-definite. Observe that the
off-diagonal entries are zero only when $\int_0^{\infty}
f_1(s)f_2(s)ds=0$. The latter integral is zero if, for instance,
we take $f_1=\varphi-\psi$ and $f_2=\varphi+\psi$, where
$\|\varphi\|=\|\psi\|<\infty$. In the third case, all the entries
of $\alpha_{\infty}$ are zero but
$(\alpha_{\infty})_{1,1}=1/\|f_1\|^2$ if $\|f_1\|<+\infty$ and
$\|f_2\|=+\infty$. The remaining case is similar by symmetry.
\end{Remark}

\begin{Remark} We shall now discuss examples of kernels of order $n$,
$n\in \mathbb{N}$, which reproducing  spaces are M\"untz spaces.
We refer to \cite{alili-wu-02} for proofs of results given below.
Take $f_i(s)=s^{\lambda_i}$, $i=1, 2, \cdots $, where
$\Lambda=\{\lambda_1,\lambda_2,\cdots\}$ is a sequence of reals
such that $\lambda_i\neq \lambda_j$ for $i\neq j$ and
$\lambda_i>-1/2$. For a fixed $n<\infty$, the kernel $k_n$ defined
by $k_n(t,s)=0$ if $s>t$ and
\begin{equation}
\label{eqn:M-K33-2}k_n(t,s)=t^{-1}\sum_{j = 1}^n a_{j,n}
(s/t)^{\lambda_j},\quad a_{j,n} = \frac{\prod_{i = 1}^n (\lambda_i
+ \lambda_j + 1)}{\prod_{i=1, i\neq j}^{n}(\lambda_i -
\lambda_j)}, \quad j=1,..., n,
\end{equation}
if $0 < s \leq t$, is a Goursat-Volterra kernel of order $n$. Its
reproducing Gaussian space, at time $t>0$, is $Span\{ \int_0^t s^i
dB_s; i=1,2,\cdots,n \}$. Going back to the Gramian matrix
 $(m_t, t\geq 0)$, observe that it has the entries
\begin{equation*}
\left( m_t \right)_{i, j}=({\lambda_i+\lambda_j
+1})^{-1}{t^{\lambda_i+\lambda_j +1}}, \qquad i, j=1,\cdots, n.
\end{equation*}
Thus if $t=1$ then $m_1$ is a Cauchy matrix. When $\lambda_i=c i$,
for some constant $c\neq 0$, and $n=\infty$, $m_1$ is the
well-known Hilbert matrix. Note that because $||f_i||= +\infty$,
$i=1,\cdots, n $, we have $\alpha_{\infty}\equiv 0$. So we have
$\phi_{i}(t)=a_{i,n}t^{-\lambda_i-1}$, $i=1,2,\cdots, n$.
Furthermore, the entries of $\alpha_t$ are given by
\begin{equation*}
(\alpha_t)_{i,j}=
a_{i,n}a_{j,n}(\lambda_i+\lambda_j+1)^{-1}t^{-\lambda_i-\lambda_j-1},
\qquad i,j=1,\cdots n,
\end{equation*}
which follows from the expression of the kernels when compared
with Theorem \ref{Theorem29}. Note that $\alpha_t$, for $t\neq 1$,
can easily be constructed from $\alpha_1$ which is known and can
be found in \cite{Schechter}. Finally, we mention that some
results are obtained about infinite order kernels in the M\"untz
case, see \cite{alili-wu-02} and \cite{hm-2004}.
\end{Remark}

\begin{Remark} \label{Hardy} Observe that we can write
\begin{equation*}
\Sigma(B)_t= \int_0^{\infty} (I-K^{*}_{f})1_{[0,t]}(u) \, dB_u,
\qquad 0<t<\infty
\end{equation*}
where $K^{*}_{f}$ is the adjoint of the bounded integral operator
$K_{f}$ defined on $L^2_{loc}(\mathbb{R}_+)$ by
\begin{equation*}
K_{f}\alpha(t)=\int_0^t k(t,r) \alpha(r) \, dr, \qquad \alpha\in
L^2_{loc}(\mathbb{R}_+).
\end{equation*}
That $I-K_{f}$ is a partial isometry, with initial subspace
$L^2_{loc}(\mathbb{R}_+)\ominus \hbox{Span}\{f\}$ and final
subspace $L^2_{loc}(\mathbb{R}_+)$, follows from the generalized
Hardy inequality
\begin{equation*}
\|K_{g}\alpha \|\leq 2\|\alpha\|, \qquad \alpha\in
L^2_{loc}(\mathbb{R}_+).
\end{equation*}
Consequently, the operator $I-K^{*}_f$, when defined on
$L^{2}_{loc}(\mathbb{R}_+)$, is isometric which implies the
statement of the Theorem \ref{hhm-thm}. For the above results, we
refer to \cite{hhm}. We also refer to the comments of Section 3
therein because here we are working with $L^2_{loc}(\mathbb{R}_+)$
instead of $L^2_{loc}([0,1])$.
\end{Remark}

\begin{Remark} Many authors work under the condition
\begin{equation}\label{integrability}
\int_0^t \left( \int_0^u k^2(u,v) \, dv \right)^{1/2} du < \infty
\end{equation}
for all $t> 0$, which is sufficient for $\Sigma(B)$, where $B$ is
a standard Brownian motion, to be well-defined, see for instance
\cite{fwy99}. However, condition (\ref{integrability}) is too
strong for $\Sigma(B)$ to be well-defined. To see that let us fix
$b\in L^2_{loc}(\mathbb{R}_+)$. The associated Goursat-Volterra
kernel of order $1$ is then found to be
\begin{equation*}
k(t,v)=b(t)b(v) \left/\int_0^{t}b^2(r) \, dr \right..
\end{equation*}
It satisfies (\ref{integrability}) if and only if $\int_0^{t}
|b(s)|/(\int_0^s b^2(r) \, dr)^{1/2} \, ds <\infty$ for all
$t<\infty$. For example, the kernel associated to
$b(t)=t^{-1}e^{-1/t}$ fails to satisfy (\ref{integrability}).
\end{Remark}


\section{{\bf On some singular linear stochastic differential
equations}}\label{sec34} Consider the singular linear stochastic
equation (\ref{ssde}). Our interest lies in the set of all its
continuous semi-martingale solutions which may be defined on a
possibly enlarged space. For a particular solution $X$, we recall
that (\ref{ssde}) is well-defined in the sense that
(\ref{condition-2}) holds. If we set $W=\Sigma(B)$, where $B$ is a
Brownian motion, then the set includes at least two solutions
which one shall now briefly describe. First, $B$ is a solution.
Second, there is a solution which is defined on $\mathbb{R}_+$ and
coincides with the $f$-generalized bridge over its life time. The
latter process, denoted by $\left(B_u^{y}, u\leq t_1\right)$, for
some $t_1>0$ and a column vector of reals $y$, is defined by
\begin{equation*}
B^{y}_u=B_u-\psi^*(u)\cdot \int_0^{t_1}f(s) \, dB_s+\psi^*(u)\cdot
y,\qquad 0<u<t_1,\end{equation*} where $\psi$ is the unique
solution to the linear system
\begin{equation*}
\int_0^u f(s) \, ds=\psi(u)\cdot \int_0^{t_1}f(s)\cdot f^*(s) \,
ds=\psi(u)\cdot m_{t_1}, \qquad 0<u<t_1.
\end{equation*}
Thus $\psi(u) = \alpha_{t_1} \cdot \int_0^{u}f(s) \, ds$ which
implies that $\int_0^{t_1} f(s) \, dB^{y}_s = y$, since
$\alpha_{t_1}$ is the inverse of $m_{t_1}$. This is the reason why
the above process is called an $f$-generalized bridge over
$[0,t_1]$ with endpoint $y$. Now, we have $\Sigma(B^{y}) =
\Sigma(B)$ which is true because $\Sigma$ is linear and
$\Sigma(\int_0^{\cdot} f(r) \, dr) \equiv 0$ since
$f(t)=\int_0^{t} k(t,v)f(v) \, dv$ for all $0 < t < \infty$. This
shows that $B^{y}$ is also a solution to (\ref{ssde}) which, in
fact, is a non-canonical decomposition. For further results on
these processes, such as their canonical decomposition in their
own filtrations, we refer to \cite{alili00}. Now, we consider
equation (\ref{ssde}) where the driving Brownian motion $W$ is
taken to be arbitrary.

\begin{Theorem}
\label{thm25} $1)$ $X$ solves equation (\ref{ssde}) if and only if
there exists a random vector $Y = (Y_1, \cdots, Y_n)^*$ such that
\begin{equation}\label{form}
X = X^0 + \int_0^{\cdot} f^*(u) \, du \cdot Y
\end{equation}
 where
\begin{equation*}
 X^0 = W -
\int_0^{\cdot} \! \int_u^{\infty} \phi^*(v) \cdot f(u) \, dW_v \,
du.
\end{equation*}
In terms of $X$, $Y$ is given by $\displaystyle Y= \lim_{t \to
\infty}
\alpha_t\cdot \int_0^t f(u) \, dX_u$. \\
$2)$ $X^0$ is a Brownian
motion if and only if $\alpha_{\infty} \equiv 0$. In case
$\alpha_{\infty} \not\equiv 0$, a process $X$ solving equation
(\ref{ssde}) is a Brownian motion if and only if  $Y$ is centered
Gaussian with covariance matrix $\alpha_{\infty}$ and is
independent ${\mathcal F}_{\infty}^{X^0}$.
\end{Theorem}

\begin{proof} 1) We proceed by checking first that
$X_t^0$ is a particular solution to (\ref{ssde}). Using the
stochastic Fubini theorem, found for instance in \cite{protter92},
we perform the decompositions
\begin{eqnarray*}
&   & X_t^0 - \int_0^t \! \int_0^u k(u,v) \, dX_v^0 \, du \\
& = & W_t - \int_0^t \! \int_u^{\infty} k(v,u) \, dW_v \, du -
\int_0^t \!
\int_0^u k(u,v) \left( dW_v - \int_v^{\infty} k(\rho,v) \, dW_{\rho} \, dv \right) du \\
& = & W_t - \int_0^t \! \int_u^{\infty} k(v,u) \, dW_v \, du - \int_0^t \! \int_0^u k(u,v) \, dW_v \, du \\
&   & \qquad + \int_0^t \! \int_0^u \! \int_0^{\rho} k(u,v)
k(\rho,v) \, dv \, dW_{\rho} \, du + \int_0^t \! \int_u^{\infty}
\! \int_0^u k(u,v) k(\rho,v) \, dv \, dW_{\rho} \, du.
\end{eqnarray*}
Since $k$ is self-reproducing, the last four terms in the last
equation cancel showing that $X_t^0$ solves (\ref{ssde}). Next, if
$X$ is a solution then by setting $X = X^0 + Z$ we see that $Z$
has to satisfy
\begin{equation*}
dZ_r = \int_0^r k(r,v) \, dZ_v \, dr, \qquad 0<r<\infty.
\end{equation*}
 Multiplying
both sides by $f(r)$ and integrating with respect to $r$, along
$[0,t]$, yields
\begin{eqnarray}\nonumber
\int_0^t f(v) \, dZ_v &=&\int_0^t f(v) \phi^*(v)\cdot
\int_0^v f(r) \, dZ_r \, dv \\
\nonumber &=& \int_0^t m_v \cdot \phi(v) \phi^*(v)\cdot
\int_0^v f(r) \, dZ_r \, dv \\
\nonumber &=&-\int_0^t m_v \cdot \frac{d}{dv} \, \alpha_v\cdot
\int_0^v f(r) \, dZ_r \, dv
\end{eqnarray}
where we used the expression of $\alpha'$ given in the proof of
Theorem \ref{Theorem29} to obtain the last equality. Because
$\alpha$ is the inverse
 of $m$, the latter relation can be written as $\frac{d}{dt}
\, \alpha_t\cdot \int_0^t f(s) \, dZ_s  =0$. This, when
integrated, yields $ \alpha_t\cdot \int_0^t f(s) \, dZ_s = Y$ for
some random vector $Y$. Hence $\int_0^t f(r) \, dZ_r = m_t\cdot Y
$ which implies that $Z_t=Y^*\cdot \int_0^{t} f(s) \, ds$. This
completes the proof of the first part of the first assertion. For
the second part, by using Theorem \ref{Theorem29} we obtain
\begin{eqnarray*}
\phi(t) \, dW_t
& = & \phi(t) \, dX_t - \phi(t) \phi^*(t)\cdot  \int_0^t f(u) \, dX_u \, dt \\
& = &  \alpha_t\cdot d \left(\int_0^t f(u) \, dX_u\right) -  \phi(t) \phi^*(t)\cdot \int_0^t f(u) \, dX_u \, dt \\
& = & d \left( \alpha_t\cdot \int_0^t f(u) \, dX_u \right).
\end{eqnarray*}
Integrating on both sides over $[s,t]$ we obtain
\begin{equation*} \int_s^t \phi(u) \, dW_u =  \alpha_t\cdot \int_0^t f(u) \, dX_u -
\alpha_s\cdot \int_0^s f(u) \, dX_u.\end{equation*} Next, observe
that as $t \to \infty$ the left hand side converges almost surely.
So the right hand side converges as well to some limit which we
denote by $\tilde{Y}$. To be more precise, setting
\begin{equation*}
\tilde{Y}= \lim_{t \to \infty}\alpha_t\cdot \int_0^t f(u) \, dX_u,
\end{equation*}
we have shown that
\begin{equation}\label{eqn:K42}
\int_t^{\infty} \phi(u) \, dW_u = \tilde{Y} - \alpha_t\cdot
\int_0^t f(u) \, dX_u,\qquad 0< t\leq  \infty.
\end{equation}
 Consequently, we have
\begin{eqnarray*}
\int_0^t \int_u^{\infty } f^*(u)\cdot \phi(v)
\, dW_v \, du &-& {\tilde{Y}}^*\cdot \int_0^t f(u) \, du\\
 & = & \int_0^t f^*(u)\cdot \alpha(u) \cdot \int_0^u f(v) \, dX_v \, du \\
& = & \int_0^t \int_0^u \phi^*(u)\cdot f(v) \, dX_v \, du.
\end{eqnarray*}
Thus, we have
\begin{eqnarray*}
 \int_0^t \! \int_u^{\infty} k(v,u) \, dW_v \, du - \tilde{Y}^* \cdot
\int_0^t f(u) \, du &=& - \int_0^t \! \int_0^u k(u,v) \, dX_v \,
du \\
&=& W_t - X_t.
\end{eqnarray*}
 Comparing with previous calculations yields $Y = \tilde{Y}$,
$\mathbb{P}_0$-almost surely. \\
2) Theorem \ref{Theorem29} implies that
\begin{equation}\label{eqn:K43} E[X_s^0 X_t^0] = s\wedge t -
\int_0^{s\wedge t} \int_0^t f^*(r)\cdot \alpha_{\infty}\cdot f(v)
\, dv \, dr.
\end{equation}
This clearly shows that $X^0$ is a Brownian motion if and only if
$\alpha_{\infty} \equiv 0$. Next, if $X$ is as prescribed  then by
virtue of (\ref{eqn:K43}), and the fact that  $\alpha_{\infty}$ is
the covariance matrix of $Y$, we have
\begin{eqnarray*}
E[X_s X_t] & = & s\wedge t -  \int_0^{s} \int_0^t f^*(u)\cdot
\alpha_{\infty}\cdot
f(v) \, dv \, du \\
&+&  \int_0^s \int_0^t E \left[ (Y^*\cdot f(u)) (Y^*\cdot f(v))
\right] \, dv \, du \\
&=& s\wedge t.
\end{eqnarray*}
Because $X$ is  a continuous Gaussian process we conclude that it
is a Brownian motion. Conversely, if $X$ is a Brownian solution to
(\ref{ssde}) then it has to be of the  form  (\ref{form}). By
virtue of the orthogonal properties of Goursat-Volterra transform,
we see that $Y$ is independent of ${\mathcal F}_t^{\Sigma(X)} =
{\mathcal F}_t^W$ for any fixed $t
> 0$. Next, by letting $t$ go to $\infty$, we get that $Y$
is independent of ${\mathcal F}_{\infty}^{X^0} \subseteq {\mathcal
F}_{\infty}^W$. Thus,  $Y$ is Gaussian vector, with covariance
matrix $\alpha_{\infty}$, which is independent of ${\mathcal
F}_{\infty}^{X^0}$ as required.
\end{proof}

Thanks to the importance of the symmetric matrix
$\alpha_{\infty}$, for instance in Theorem \ref{Theorem29}, it is
natural to look for a description of its structure. The following
result, which is hidden in the proof of Theorem \ref{thm25}, gives
a necessary and sufficient condition for a column or a row to be
zero.
\begin{Corollary} \label{description} For $1 \le i \le n$, $(\alpha_{\infty})_{i,j} = (\alpha_{\infty})_{j,i} =
0$ for all $j$, if and only if $\| f_i \| = \infty$.
\end{Corollary}

\begin{proof}
For a fixed $t>0$, $\alpha_t$ is the covariance matrix of
$\alpha_t \cdot \int_0^t f(s) \, dB_s$. Furthermore, due to
Theorem \ref{thm25}, we conclude that $\alpha_t \cdot \int_0^t
f(s) \, dB_s$ converges to a Gaussian vector $Y$, possibly with
some null components, such that $E(Y \cdot Y^*)=\alpha_{\infty}$.
Thus, $Y_i\equiv 0$ for some $i$ if and only if
$(\alpha_{\infty})_{i,i}=0$ and if and only if $\|f_i\|=\infty$.
Now, $(\alpha_{\infty})_{i,i}=0$ if and only if
$(\alpha_{\infty})_{i,j} = 0$ for all $j$. In order to see that,
we let $t\rightarrow \infty$ and use continuity in the well-known
inequality $|(\alpha_{t})_{i,j}|^2 \leq (\alpha_{t})_{i,i}
(\alpha_{t})_{j,j}$ valid for symmetric positive definite
matrices.
\end{proof}

Now, we take a look at the orthogonal decompositions of
filtrations which arise from Goursat-Volterra transforms and
provide their interpretation.

\begin{Corollary}\label{decomposition}
The orthogonal decomposition given by (\ref{eqn:C70}) holds true.
Furthermore, the progressive decomposition
\begin{equation*}
{\mathcal F}_t^{B} = {\mathcal F}_t^{\Sigma(B)} \otimes  \sigma
\left( Y- \int_t^{\infty} \phi(u) \, d\Sigma(B)_u \right), \qquad
0<t<\infty
\end{equation*}
 holds true, where $Y\equiv 0$ if $\alpha_{\infty}\equiv 0$ and
$Y$ is a Gaussian vector independent of ${\mathcal
F}_{\infty}^{\Sigma(B)}$ with covariance matrix $\alpha_{\infty}$
otherwise. Thus,  we have ${\mathcal F}_{\infty}^B = {\mathcal
F}_{\infty}^{\Sigma(B)}$ in case $\alpha_{\infty}\equiv 0$ and
${\mathcal F}_{\infty}^B = {\mathcal
F}_{\infty}^{\Sigma(B)}\vee\sigma\{ Y\}$ otherwise.
\end{Corollary}

\begin{proof} For a fixed $t> 0$, Theorem \ref{thm25} implies that
\begin{equation*} B_t = \Sigma(B)_t - \int_0^t \! \int_u^{\infty} k(v,u) \,
d\Sigma(B)_v \, du + Y^*\cdot \int_0^t f(u) \, du \end{equation*}
where $Y$ is a Gaussian vector with covariance $\alpha_{\infty}$
which is independent of ${\mathcal F}_{\infty}^{\Sigma(B)} $.
Hence, we have
\begin{equation*}
\int_0^t f(u) \, dB_u = m_t\cdot \left( Y - \int_t^{\infty}
\phi(u) \, d\Sigma(B)_u \right)
\end{equation*}
 which gives
\begin{equation*}
 \sigma \left\{\int_0^t f(u)
\, dB_u \right\} =\sigma \left\{ Y - \int_t^{\infty} \phi(u) \,
d\Sigma(B)_u \right\}.
\end{equation*}
 This implies the first assertion while the last one follows
by letting $t$ tend to $+\infty$.
\end{proof}

\begin{Remark} Recall that $\mathcal{F}_0^{B}$ and $\mathcal{F}_0^{\Sigma(B)}$ are trivial.
So by letting $t$ converge to
 $0$, in Corollary \ref{decomposition},
we see that
 $\phi^* \in L^2([\varepsilon, \infty)^n)$ for all $\varepsilon
> 0$ but $\phi_i \not\in L^2((0,+\infty))$, for $i=1, \cdots n$.
This fact can also be shown by a combination of Theorem
\ref{Theorem29} and the inequality
\begin{equation*}
(\alpha_{t})_{i,i} \geq {1}/{(m_t)_{i,i}}={1}/{\|f_i\|^2}
 \end{equation*}
which follows from the orthogonal diagonalization of $m_t$ and may
be found in Exercise 8, p. 274, in \cite{D.Harville}.
\end{Remark}

\begin{Remark} It is clear that if the choice of the vector $f$ allows the use of
integration by parts for the integrand in the right hand side of
(\ref{ssde}) then we obtain a stochastic differential equation
which does not involve a stochastic integral. For instance, that
is the case for the examples given by P. L\'evy,  found in
\cite{Levy-56} and \cite{Levy-57}. These go back to around the
middle of last century  when stochastic integration was not yet
world-widely developed.
\end{Remark}

\section{\bf Connections to some positive martingales}
 Let $(k(t,s), t\geq s>0)$ be a Goursat-Volterra kernel of order
$n$, where $n$ is a natural number.  Assume that $f$ is a
reproducing basis for $k$, or for the associated Volterra
transform $\Sigma$, and let us keep the notations used in the
Introduction. Consider the singular stochastic differential
equation (\ref{ssde}) associated to $k$ and driven by a given
standard Brownian motion $W$. Our aim here is to describe the set
\begin{eqnarray*}
\Upsilon^{(k)}& = &\{\mathbb{P} \hbox{ is the probability law of a
continuous semi-martingale } X
\\
& & \hbox{on }\; (\mathcal{C}([0,\infty), \mathbb{R}), {\mathcal
F}_{\infty}^*) \; \hbox{solving}\; (\ref{ssde})\;
\hbox{s.t.}\; \Sigma{(X)} \hbox{ is a Brownian motion} \\
& & \hbox{and}\;  \mathcal{F}_t^{\Sigma(X)}\; \hbox{is independent
of} \int_0^{t} f(s) \, dX_s, \hbox{ for all } 0<t<\infty \}.
\end{eqnarray*}
We read from Corollary \ref{description} that
$\alpha_{\infty}\equiv 0$ if and only $||f_i||=\infty$ for all
$i$.  Now,  we are ready to state the following unified
characterization of the set $\Upsilon^{(k)}$.
\begin{Theorem}\label{thm-harmonic} If $\alpha_{\infty}\equiv 0$ then the following
assertions are equivalent
\begin{itemize}
\item[{(1)}] $\mathbb{P} \in \Upsilon^{(k)}$.

\item[{(2)}] $\mathbb{P}$ is the law of ${\displaystyle B +
Y^*\cdot \int_0^{.} f(s) ds}$, where $B$ is a standard Brownian
motion and $Y$ is a vector of random variables which is
independent of ${\mathcal F}_{\infty}^B$.

\item[{(3)}] There exists a  positive  function $h\in
C^{1,2}\left( \real_+ \times \real^n, \real_+\right)$ such that $
h(., \int_0^{.} f^*(s) \, dB_s)$ is a continuous
$(\mathbb{P}_0,{\mathcal F})$-martingale with expectation $1$,
 and
$\mathbb{P} = \mathbb{P}_0^h$ with
\begin{equation*}
 \mathbb{P}_0^h \left.
\right|_{{\mathcal F}_t} =  h \left( t, \int_0^t f^*(s) \, dB_s
\right)
 \cdot \mathbb{P}_0 \left. \right|_{{\mathcal F}_t}, \qquad 0<t<\infty,
\end{equation*}
where $\mathbb{P}_0$ stands for the Wiener measure,
\end{itemize}
\end{Theorem}

\begin{proof} We split the proof into several steps where we show
that $(1)\Longleftrightarrow (2)$ and $(2)\Longleftrightarrow
(3)$. Let us show that $(1) \Longrightarrow (2)$. Let $ \mathbb{P}
\in \Upsilon^{(k)}$.  Theorem \ref{thm25} implies that there exist
a vector $Y$ such that $\mathbb{P}$ is the law of $X_t^0 +
Y^*\cdot \int_0^t f(u) \, du$. That combined with the assumption
$\alpha_{\infty}\equiv 0$ leads to the fact that  $X^0$ is a
Brownian motion. Hence, it suffices to show that $Y$ is
independent of $X^0$. From (\ref{eqn:K42}) we see that $Y=
\int_t^{\infty} \phi(u) \, dB_u + \alpha_t\cdot  \int_0^t f(u) \,
dX_u$,
 the vector $\int_0^t f(u) dX_u$ is independent of $B$ and,
consequently, it is also independent of $X^0$. Thus, whenever
$Z\in L^2({\mathcal F}_{\infty}^{X^0})$, for any fixed $t\geq 0$,
we have
\begin{eqnarray*}
E \left[ E \left[ Z \left| {\mathcal F}_t^{X^0} \right. \right]
\phi \left( Y - \int_t^{\infty} \phi(u)
 dB_u  \right) \right]=E \left[ Z \right] E \left[ \phi \left( Y- \int_t^{\infty}
\phi(u) \, dB_u\right) \right]
\end{eqnarray*}
for any  bounded function $\phi: \real^n \to \real$. By letting $t
\rightarrow \infty$ we conclude that  $E \left[ Z \cdot \phi(Y)
\right] = E[Z] E \left[ \phi(Y) \right]$ which implies the
required independence. We shall now show that $ (2)
\Longrightarrow (1)$. To this end, let $k$ be a Goursat-Volterra
kernel.  Denote by $f$  a reproducing basis associated to $k$ and
 put $X_t= B_t + Y^*\cdot \int_0^{t} f(s) \, ds$ for $t>0$.
For a fixed $t>0$, because $\int_0^{t} f(u) \, dB_u \in
\Gamma_{t}^{(k)}$,  we can write
\begin{eqnarray*}
X_t - \int_0^{t} \! \int_0^u k(u,v) \, dX_v \, du &=& B_t -
\int_0^{t} \! \int_0^u k(u,v) \, dB_v \, du = \Sigma(B)_t
\end{eqnarray*}
which, of course, a Brownian motion. Furthermore, using once more
the above argument we can easily see that $ \int_0^t f(u) \, dX_u$
is independent of $\mathcal{F}_t^{\Sigma(B)}$.  Next, we deal with
$(2) \Longrightarrow (3)$. Denote by $\nu(dy)$ the distribution of
$Y$. For any measurable functional $\phi$, we then have that
\begin{eqnarray*}
&   & E \left[ \phi \left( B_s +  Y^*\cdot \int_0^s f(u) \, du:
s \le t \right) \right] \\
& = & \int_{\real^n} E \left[ \phi \left( B_s +  y^*\cdot \int_0^s
f(u) \, du: s \le t \right) \right] \nu(dy) \\
& = & \int_{\real^n} E \left[ \exp \left(\int_0^t  y^*\cdot f(u) \,
dB_u - \frac12 \int_0^t \left(  y^*\cdot f(u) \right)^2 \, du
\right) \phi(B_s: s \le t) \right] \nu(dy)
\end{eqnarray*}
where the last equality is obtained by Girsanov theorem. The
required space-time harmonic function is thus given on
$\mathbb{R}^+\times \mathbb{R}^n$ by
\begin{equation*} h(t, x) = \int_{\real^N} \exp \left(  y^* \cdot x - \frac12
\int_0^t \left( y^*\cdot f(s) \right)^2 \, ds \right) \nu(dy).
\end{equation*}
It remains to show that  $(3) \Longrightarrow (2)$. For fixed
$0<u\leq t<+\infty$,  set $\psi(u,t) = \alpha_t\cdot \int_0^u f(s)
\, ds$. Let us write the obvious decomposition
\begin{eqnarray*}
B_u &=& \left( B_u - \psi^*(u,t)\cdot \int_0^t f(s) \, dB_s
\right) +  \psi^*(u,t)\cdot \int_0^t f(s) \, dB_s
\end{eqnarray*}
and denote by $H_u^t$ the first term of its right hand side. We
observe that the process $(H_u^t, u<t)$ has has then the same law
under $\mathbb{P}_0$ as under $\mathbb{P}_0^h$. Next, to simplify
notations, write
\begin{equation*}
\hat{H}_u^t= \psi^*(u,t)\cdot \int_0^t f(s) \, dB_s
 =  Y^*_t\cdot \int_0^u f(s) \, ds,
\end{equation*}
where we set $Y^*_t=\alpha_t \cdot\int_0^t f(r)dB_r$. For any $ 0
\le s \le u \le t$, we have $E[H_s^t H_u^t] = s - \psi^*(s,t)\cdot
\int_0^u f(v) \, dv$ and $ \psi^*(s,t)\cdot \int_0^u f(v) \, dv =
\int_0^u f^*(v) \, dv \cdot \alpha_t \cdot \int_0^s f(r) \,
d\rightarrow 0$ as $t\rightarrow \infty$ because
$\alpha_{\infty}\equiv 0$. We conclude that the convergence in
distribution $H_{\cdot}^t\rightarrow B^{(h)}_{\cdot}$ holds, where
$B^{(h)}$ is a $\mathbb{P}_0^h$-Brownian motion. That implies the
convergence of $\hat{H}_{.}^t$ as well to a finite limit. But that
can happen if and only if $Y^*_t$ converges to a finite limit
which we denote by $Y^*$. Finally, from the above arguments we see
that $Y^*$ is independent of $\mathcal{F}_{\infty}^{B^{(h)}}$
which ends the proof.
\end{proof}

\begin{Remark}
Unfortunately, for the case $\alpha_{\infty} \not\equiv 0$, the
second statement in the above theorem is too strong. For example,
$X^0$ satisfies the assertion (1) but it is easily seen  that it
does not satisfy (2). The implications  $(2) \Longrightarrow (1)$
and $(2) \Longrightarrow (3)$ still work in this case. We also can
replace $B$ by $X^0$ in statement (2) and prove  that $(1)
\Longleftrightarrow (2)$ still holds true. However, $(2)
\Longrightarrow (3)$ fails.
\end{Remark}


\noindent{\bf Acknowledgment:} This work was partly supported by
the Austrian Science Foundation (FWF) under grant
Wittgenstein-Prize Z36-MAT. We would like to thank W.
Schachermayer and his team for their warm reception. The second
author is greatly indebted to the National Science Council Taiwan
for the research grant NSC 96-2115-M-009-005-MY2. We thank  Th.
Jeulin, R. Mansuy and M. Yor for fruitful discussions. Finally, we
are grateful for anonymous referees for several reports on
previous versions of this paper which lead to its improvement.


\vspace{5mm} \noindent{\footnotesize $^{(1)}$ Department of
Statistics, University of Warwick, CV4 7AL, Coventry, UK. L.alili@warwick.ac.uk\\
$^{(2)}$ Department of Applied Mathematics,
             National Chiao Tung University,
             No. 1001, Ta-Hsueh Road,
             300 Hsinchu, Taiwan. ctwu@math.nctu.edu.tw}

\end{document}